\documentclass{article}
\begin{document}
\newtheorem{proposition}{Proposition}[section]
\newtheorem{definition}{Definition}[section]
\newtheorem{lemma}{Lemma}[section]
\newcommand{\xl}{\stackrel{\rightharpoonup}{\cdot}}
\newcommand{\xr}{\stackrel{\leftharpoonup}{\cdot}}
\newcommand{\xlplus}{\stackrel{\rightharpoonup}{+}}
\newcommand{\xrplus}{\stackrel{\leftharpoonup}{+}}
\newcommand{\xluplus}{\stackrel{\rightharpoonup}{\uplus}}
\newcommand{\xruplus}{\stackrel{\leftharpoonup}{\uplus}}
\newcommand{\xlodot}{\stackrel{\rightharpoonup}{\odot}}
\newcommand{\xrodot}{\stackrel{\leftharpoonup}{\odot}}
\newcommand{\Ll}{\stackrel{\rightharpoonup}{L}}
\newcommand{\Lr}{\stackrel{\leftharpoonup}{L}}
\newcommand{\Rl}{\stackrel{\rightharpoonup}{R}}
\newcommand{\Rr}{\stackrel{\leftharpoonup}{R}}

\title{\bf How to Expand the Zariski Topology}
\author{Keqin Liu\\Department of Mathematics\\The University of British Columbia\\Vancouver, BC\\
Canada, V6T 1Z2}
\maketitle

\begin{abstract} We introduce the notion of a Hu-Liu prime ideal in the context of left commutative rngs, and establish the contravariant functor from the category of left commutative rngs into the category of topological spaces.\end{abstract}

It is well known that new points must be introduced in order to expand algebraic geometry over algebraically closed fields into Grothendieck's scheme theory over commutative rings. We believe that the idea of adding new points to an old space is still essential for the attempt to expand algebraic geometry over algebraically closed fields into a kind of geometry over a class of 
non-commutative rings. Clearly, whether we can use the natural idea successfully depends on whether we can find new points satisfactorily.  Since points in Grothendieck's scheme theory are prime ideals, the problem of finding satisfactory new points is how to choose some classes of rings to get a satisfactory generalization of prime ideals. The purpose of this paper is to give a solution to the problem. Our solution is based on the notion of the additive halo introduced in \cite{Liu4}, which comes from some facts obtained in our attempt to generalize the Lie correspondence between connected linear Lie groups and linear Lie algebras. 

\medskip
After choosing a class of rings called left commutative rngs, we introduce the notion of a Hu-Liu prime ideal, characterize the nil radical by using Hu-Liu prime ideals, and establish the contravariant functor from the category of left commutative rngs into the category of topological spaces.

\bigskip
\section{Definitions}

Following \cite{J}, the word ``rng'' means a ring which is not assumed to have an identity, and the word ``ring'' always means a ring with an identity. 

\medskip
We begin this section with the definition of a left commutative rng.

\begin{definition}\label{def1.1} A rng $(\, R, \, + , \, \cdot \,)$ is called a {\bf left commutative rng} if $R$ satisfies the following four conditions.
\begin{description}
\item[(i)] The associative product $\cdot$ is {\bf left commutative}; that is, 
\begin{equation}\label{eq1}
xyz=yxz \quad\mbox{for $x$, $y$, $z\in R$.}
\end{equation}
\item[(ii)] There exists a {\bf left identity} $1^{\ell}$ of $R$ such that
\begin{equation}\label{eq2}
1^{\ell}x=x \quad\mbox{for $x\in R$.}
\end{equation}
\item[(iii)] There exists a binary operation $\,\sharp\,$ called the {\bf local product} on the {\bf additive halo} 
$$\hbar ^+(R): =\{\, x \,|\, \mbox{$x\in R$ and $x1^{\ell}=0$} \,\}$$ 
such that $(\, \hbar ^+(R) , \, + , \, \,\sharp\, \, )$ is a commutative ring with an identity $1^\sharp$. The identity $1^\sharp$ of the ring $\hbar ^+(R)$ is called the {\bf local identity} of $R$.
\item[(iv)] The two associative binary operations $\cdot$ and $\,\sharp\,$ satisfy the {\bf Hu-Liu  triassociative law}:
\begin{equation}\label{eq3}
(x\alpha ) \,\sharp\, \beta=x(\alpha \,\sharp\, \beta),
\end{equation}
where $\alpha$, $\beta\in \hbar ^+(R)$ and $x\in R$.
\end{description}
\end{definition}

Since the local product $\,\sharp\,$ is commutative, (\ref{eq3}) is equivalent to
\begin{equation}\label{eq4}
\alpha  \,\sharp\, (x\beta)=x(\alpha \,\sharp\, \beta) \quad\mbox{for $\alpha$, 
$\beta\in \hbar ^+(R)$ and $x\in R$.}
\end{equation}
The equations (\ref{eq3}) and (\ref{eq4}) consist of a version of the Hu-Liu  triassociative law
introduced in \cite{Liu5}.

\medskip
A left commutative rng is sometimes denoted by $(\, R, \, + , \, \cdot , \, \,\sharp\, \,)$. If $R$ is a left commutative rng, then the {\bf bar-unit set} 
$\hbar ^\times (R)$ of  $R$ is defined by
$$\hbar ^\times (R): =\{\, b \,|\, \mbox{ $bx=x$ for $x\in R$} \,\}.$$ 

\begin{definition}\label{def1.2} Let $R$ and $\bar{R}$ be two left commutative rngs with the local identity $1^\sharp$ and $\bar{1}^\sharp$, respectively. A map $\phi: R\to \bar{R}$ is called a {\bf left commutative rng homomorphism} if 
\begin{eqnarray}
\label{eq5} \phi (x+y)&=&\phi (x)+\phi (y),\\
\label{eq6} \phi (xy)&=&\phi (x)\phi (y),\\
\label{eq7} \phi \left(\hbar ^{\times} (R)\right)\cap \hbar ^{\times} (\bar{R})&\ne& \emptyset,\\
\label{eq8} \phi (\alpha \,\sharp\, \beta )&=&\phi (\alpha) \,\sharp\, \phi (\beta),\\
\label{eq9} \phi (1^\sharp)&=&\bar{1}^\sharp,
\end{eqnarray}
where $x$, $y\in R$, and $\alpha$, $\beta\in \hbar ^+(R)$. A bijective left commutative rng homomorphism is called a {\bf left commutative rng isomorphism}.
\end{definition}

Note that (\ref{eq8}) is well-defined because of the following fact:
$$
\mbox{(\ref{eq5}), (\ref{eq6}) and (\ref{eq7}) $\Rightarrow$ 
$\phi (\hbar ^+(R))\subseteq \hbar ^+(\bar{R})$.}
$$

\medskip
If $(\, R, \, + , \, \cdot , \, \,\sharp\, \,)$ is a left commutative rng with a left identity $1^{\ell}$, then the left identity $1^{\ell}$ induces a decomposition of $R$:
\begin{equation}\label{eq10}
R=R_0 \oplus R_1 \quad\mbox{(as Abelian groups)},
\end{equation}
where
$$
R_0 : =R1^{\ell}= \{\, x1^{\ell} \,|\, x\in R \, \} \quad\mbox{and}\quad
R_1 : = \hbar ^+(R) = \{\, x \,|\, \mbox{$x\in R$ and $x1^{\ell}=0$} \,\}.
$$
$R_0$ and $R_1$ are called the {\bf even part} of $R$ induced by the left identity $1^{\ell}$ and the {\bf odd part} of $R$, respectively. If $x\in R$, then
$$
x=x_0 + x_1, \qquad\mbox{$x_0\in R_0$ and $x_1\in R_1$}
$$
by (\ref{eq10}). $x_0$ and $x_1$ are called the {\bf even component} and  the {\bf odd component} of $x$ induced by the left identity $1^{\ell}$, respectively. We also say that $x_\varepsilon$ is the {\bf $\varepsilon$-component} of $x$ induced by the left identity $1^{\ell}$, where $\varepsilon =0$, $1$.

\medskip
Let $(\,R , \, + , \, \cdot , \,\sharp\, \,)$ be a left commutative rng. A subgroup $I$ of the additive group $(\,R , \, +  \,)$ is called an {\bf ideal} if $RI\subseteq I$, $IR\subseteq I$ and $I\cap \hbar ^+(R)$ is an ideal of the ring $(\, \hbar ^+(R) , \, + , \, \,\sharp\, \,)$.
If $\hbar ^+(R)\ne 0$, then every left commutative rng $R$ always has three distinct ideals: $0$, $\hbar ^+(R)$ and $R$. An ideal $I$ of $R$ respects to the decomposition (\ref{eq10}); that is,
$$
I=I_0\oplus I_1,
$$
where $I_0=I\cap R1^{\ell}$ and $I_1=I\cap \hbar ^+(R)$ are called the {\bf even part} and {\bf odd part} of $I$ induced by the left identity $1^{\ell}$.

Let $I$ be an ideal of a triring $(\, R, \, + , \, \xl , \, \xr , \,\sharp\, \,)$. We define a binary operation $\cdot$ on the quotient group
$$\displaystyle\frac{R}{I}:=\left\{\, x+I \, |\, x\in I \, \right\}$$
by
\begin{equation}\label{eq11}
(x+I)\cdot (y+I):=xy +I,
\end{equation}
where $x$, $y\in  R$. The well-defined binary operation above makes the quotient group $\displaystyle\frac{R}{I}$ into a rng with a left identity $1^{\ell}+I$, where $1^{\ell}$ is a left identity of $R$. The additive halo of $\displaystyle\frac{R}{I}$ is given by
\begin{equation}\label{eq12}
\hbar ^+\left(\displaystyle\frac{R}{I}\right)=
\displaystyle\frac{\hbar ^+(R) +I}{I}=\left\{\, \alpha +I \,|\, \alpha\in \hbar ^+(R)\,\right\}.
\end{equation}

We now define a local product on $\hbar ^+\left(\displaystyle\frac{R}{I}\right)$ by
\begin{equation}\label{eq13}
(\alpha +I) \,\sharp\, (\beta +I) : =\alpha \,\sharp\, \beta +I,
\end{equation}
where $\alpha$, $\beta\in \hbar ^+ (R)$.
One can check that (\ref{eq13}) is well-defined, the two binary operations defined by 
(\ref{eq11}) and (\ref{eq13}) satisfy Hu-Liu  triassociative law, and 
$\left( \, \hbar ^+\left(\displaystyle\frac{R}{I}\right), \, + , \, \,\sharp\, \,\right)$ is a ring with the identity $1^\sharp +I$, where $1^\sharp$ is the local identity of $R$. Therefore, the quotient group $\displaystyle\frac{R}{I}$ becomes a  left commutative rng under (\ref{eq11}) and (\ref{eq13}), which is called the {\bf quotient left commutative rng} of $R$ with respect to the ideal $I$.

\medskip
The following definition gives the counterpart of a prime ideal in the context of left commutative rngs.

\begin{definition}\label{def1.3} Let $(\, R, \, + , \, \cdot , \, \,\sharp\, \,)$ be a left commutative rng. An ideal $P$ of $R$ is called a {\bf Hu-Liu prime ideal} if $P\ne R$ and the following three  conditions are satisfied.
\begin{description}
\item[(i)] For $x$, $y\in R$, we have 
\begin{equation}\label{eq14}
xy\in P+ \hbar ^+(R) \,\Rightarrow \,\mbox{$x\in P+ \hbar ^+(R)$ or $y\in P+ \hbar ^+(R)$.}
\end{equation}
\item[(ii)] For $x$, $y\in R$, we have 
\begin{equation}\label{eq15}
xy\in P \,\Rightarrow \,\mbox{$x\in P+ \hbar ^+(R)$ or $y\in P$.}
\end{equation}
\item[(iii)] $P\supseteq \hbar ^+(R)$ or $P\cap \hbar ^+(R)$ is a prime ideal of the commutative ring 
\linebreak
$(\, \hbar ^+(R) , \, + , \, \,\sharp\, \, )$.
\end{description}
\end{definition}

Let $R$ be a left commutative rng. The set of all Hu-Liu prime ideals of $R$ is called the {\bf spectrum} of $R$ and is denoted by $spec^{\sharp}R$. It is clear that
$$spec^{\sharp}R=spec^{\sharp}_0R\cup spec^{\sharp}_1R\quad\mbox{and}\quad 
spec^{\sharp}_0R\cap spec^{\sharp}_1R=\emptyset,$$
where
$$spec^{\sharp}_0R:=\{\, P \, |\, \mbox{$P\in spec^{\sharp}R$ and $P\supseteq \hbar ^+(R)$}\,\}$$
is called the {\bf even spectrum} of $R$ and 
$$spec^{\sharp}_1R:=\{\, P \, |\, \mbox{$P\in spec^{\sharp}R$ and $P\not\supseteq \hbar ^+(R)$}\,\}$$
is called the {\bf odd spectrum} of $R$.

\medskip
The next proposition gives some equivalent forms of the conditions in Definition~\ref{def1.3}.

\begin{proposition}\label{pr1.1} Let $P=P_0\oplus P_1\ne R$ be an ideal of a left commutative rng $R$ with a left identity $1^{\ell}$, where $P_{\varepsilon}=P\cap R_{\varepsilon}$, $R_0:=R1^{\ell}$ and  $R_1:=\hbar ^+(R)$.
\begin{description}
\item[(i)] $P$ is a Hu-Liu prime ideal of $R$ if and only if 
\begin{equation}\label{eq16} 
\mbox{$x_0y_{\varepsilon}\in P  \,\Rightarrow \, $ $x_0\in P_0$ or $y_{\varepsilon}\in P_{\varepsilon}$}
\end{equation}
and
\begin{equation}\label{eq17} 
\mbox{$x_1\,\sharp\, y_1\in P_1\,\Rightarrow \,$  $x_1\in P_1$ or $y_1\in P_1$,}
\end{equation}
where $x_{\varepsilon}$, $y_{\varepsilon}\in R_{\varepsilon}$ and $\varepsilon=0$, $1$.
\item[(ii)] If $P\supseteq \hbar ^+(R)$, then $P$ is a Hu-Liu prime ideal of the left commutative rng $R$ if and only if $\displaystyle\frac{P}{\hbar ^+(R)}$ is a prime ideal of the commutative ring $\displaystyle\frac{R}{\hbar ^+(R)}$.
\item[(iii)] (\ref{eq14}) holds if and only if $P_0$ is a prime ideal of the commutative ring $R_0$.
\item[(iv)] (\ref{eq15}) holds if and only if the left $\displaystyle\frac{R_0}{P_0}$-module $\displaystyle\frac{R}{P}$ is faithful, where the left $\displaystyle\frac{R_0}{P_0}$-module action is defined by 
$$
(x_0+P_0)(y+P):= x_0y+P \quad\mbox{$x_0\in P_0$ and $y\in R$.}
$$
\item[(v)] If $P\not\supseteq \hbar ^+(R)$, then $P_1=P\cap \hbar ^+(R)$ is a prime ideal of the commutative ring $(\, \hbar ^+(R) , \, + , \, \,\sharp\, \, )$ if and only if 
$\left(\, \hbar ^+\left(\displaystyle\frac{R}{P}\right) , \, + , \, \,\sharp\, \, \right)$ is a domain.
\end{description}
\end{proposition}

\medskip
\noindent
{\bf Proof} They are direct consequences of Definition~\ref{def1.3}.

\hfill\raisebox{1mm}{\framebox[2mm]{}}

\bigskip
\section{Nil Radicals}

\medskip
Let $(\, R, \, + , \, \cdot , \, \,\sharp\, \,)$ be a left commutative rng with a left identity $1^{\ell}$ and a local identity $1^\sharp$. 

For $a\in R$ and $\alpha \in \hbar ^+(R)$, we define the {\bf $n$th power} $a^n$ and the 
{\bf $n$th local power} $\alpha ^{\sharp n}$ as follows:
$$
a^n :=\left\{ \begin{array}{ll}
1^{\ell}, & \qquad\mbox{if $n=0$;}\\
\underbrace{aa\cdots a}_n , & \qquad\mbox{if $n$ is a positive integer} \end{array}\right.
$$
and
$$
\alpha^{\sharp n} :=\left\{ \begin{array}{ll}
1^\sharp, & \qquad\mbox{if $n=0$;}\\
\underbrace{\alpha \,\sharp\, \alpha  \,\sharp\, \cdots  \,\sharp\, \alpha}_n , & \qquad\mbox{if $n$ is a positive integer.} \end{array}\right.
$$
A product $a^m(\alpha^{\sharp n})$ will be denoted by $a^m\alpha^{\sharp n}$.

\medskip
\begin{proposition}\label{pr2.1} Let $x$ be an element of a left commutative rng $(\, R, \, + , \, \cdot , \, \,\sharp\, \,)$. Let $1^{\ell}$ and $\overline{1^{\ell}}$ be two left identities of $R$. If $p_1(x)$ and $\overline{p_1}(x)$ are odd components of $x$ induced by $1^{\ell}$ and $\overline{1^{\ell}}$ respectively, then the following are equivalent:
\begin{description}
\item[(i)] $x^m=0$ and $p_1(x)^{\,\sharp\, n}=0$ for some $m$, $n\in \mathcal{Z}_{>0}$.
\item[(ii)] $x^m=0$ and $\overline{p_1}(x)^{\,\sharp\, k}=0$ for some $m$, $k\in \mathcal{Z}_{>0}$.
\end{description}
\end{proposition}

\medskip
\noindent
{\bf Proof} Since $p_1(x)=x-x1^{\ell}$ and $\overline{p_1}(x)=x-x\overline{1^{\ell}}$, we have
$\overline{p_1}(x)=p_1(x)+x(1^{\ell} -\overline{1^{\ell}})$. It follows from (\ref{eq3}) that
\begin{eqnarray}
\overline{p_1}(x)^{\,\sharp\, k}&=& \left( p_1(x)+x(1^{\ell} -\overline{1^{\ell}})\right)^{\,\sharp\, k}\nonumber\\
&=& \displaystyle\sum _{i=0}^{k} {k \choose i} p_1(x)^{\,\sharp\, i} \,\sharp\, 
\left( x(1^{\ell} -\overline{1^{\ell}})\right)^{\,\sharp\, (k-i)}\nonumber\\
\label{eq18}&=&\displaystyle\sum _{i=0}^{k} {k \choose i} p_1(x)^{\,\sharp\, i} \,\sharp\, 
\left( x^{k-i}(1^{\ell} -\overline{1^{\ell}})\right)^{\,\sharp\, (k-i)}.
\end{eqnarray}

If (i) is true, then (ii) is true for $s=m$ and $k\ge m+n-1$ by (\ref{eq18}).

\medskip
Similarly, (ii) implies (i).

\hfill\raisebox{1mm}{\framebox[2mm]{}}

\bigskip
\begin{definition}\label{def2.1} Let $(\, R, \, + , \, \cdot , \, \,\sharp\, \,)$ be a left commutative rng with a left identity $1^{\ell}$. An element $x$ of $R$ is called a {\bf nilpotent} element if 
$$
\mbox{$x^m=0$ and $x_1^{\,\sharp\, n}=0$ for some $m$, $n\in\mathcal{Z}_{>0}$,}
$$
where $x_1:=p_1(x)$ is the old component of $x$ induced by $1^{\ell}$.
\end{definition}

It is clear that Definition~\ref{def2.1} is independent of the choice of a left identity by Proposition~\ref{pr2.1}.

\begin{definition}\label{def2.2} The {\bf nil radical} of a left commutative rng 
$(\, R, \, + , \, \cdot , \, \,\sharp\, \,)$ is the set of nilpotent elements of $R$. We shall use $nilrad ^{\sharp}(R)$ or $\sqrt[\sharp]{0}$ to denote the nil radical of the left commutative rng $R$.
\end{definition}

The nil radical of a left commutative rng $R$ can be expressed as 
\begin{equation}\label{eq19} 
nilrad ^{\sharp}(R)=nilrad (R1^{\ell})  \oplus nilrad (\hbar ^+(R)),
\end{equation}
where $1^{\ell}$ is a left identity of $R$, $nilrad (R1^{\ell})$ is the nil radical of the commutative ring $(\, R1^{\ell}, \, + , \, \cdot \,)$, and $nilrad (\hbar ^+(R))$ is the nil radical of the commutative ring $(\, \hbar ^+(R), \, + , \, \,\sharp\, \,)$.

\begin{proposition}\label{pr2.2} Let $(\, R, \, + , \, \cdot , \, \,\sharp\, \,)$ be a left commutative rng.
\begin{description}
\item[(i)] The nil radical $nilrad ^{\sharp}(R)$ is an ideal of $R$.
\item[(ii)] $nilrad ^{\sharp}\left( \displaystyle\frac{R}{nilrad ^{\sharp}(R)}\right)=0.$
\end{description}
\end{proposition}

\medskip
\noindent
{\bf Proof} (i) By (\ref{eq19}), it is enough to prove that 
\begin{equation}\label{eq20} 
R\left(nilrad (\hbar ^+(R))\right)\subseteq nilrad (\hbar ^+(R)).
\end{equation}

If $x\in R$ and $\alpha \in nilrad (\hbar ^+(R))$, then $\alpha^{\sharp m}=0$ for some $m\in \mathcal{Z}_{>0}$. It follows that
$$(x\alpha)^{\sharp m}=x^m \alpha^{\sharp m}=x^m 0=0,$$
which implies (\ref{eq20}).

\medskip
(ii) If $x+nilrad ^{\sharp}(R)\in nilrad ^{\sharp}\left(\displaystyle\frac{R}{nilrad ^{\sharp}(R)}\right)$, then
\begin{equation}\label{eq21} 
\left(x+nilrad ^{\sharp}(R)\right)^m=nilrad ^{\sharp}(R)
\end{equation}
and
\begin{equation}\label{eq22} 
\left((x+nilrad ^{\sharp}(R))_1\right)^{\,\sharp\, n}=nilrad ^{\sharp}(R)
\end{equation}
for some positive integers $m$ and $n$.

By (\ref{eq21}), we get
\begin{eqnarray}
&&x^m+nilrad ^{\sharp}(R)=nilrad ^{\sharp}(R)\nonumber\\
&\Rightarrow& x^m\in nilrad ^{\sharp}(R)\nonumber\\
\label{eq23}&\Rightarrow& x^{mu}=(x^m)^u=0 \quad\mbox{for some $u\in \mathcal{Z}_{>0}$.}
\end{eqnarray}

Since $(x+nilrad ^{\sharp}(R))_1=x_1+nilrad ^{\sharp}(R)$, we get from (\ref{eq22}) that
\begin{eqnarray}
&&x_1^{\,\sharp\, n}+nilrad ^{\sharp}(R)=nilrad ^{\sharp}(R)\nonumber\\
&\Rightarrow& x_1^{\,\sharp\, n}\in nilrad ^{\sharp}(R)\nonumber\\
\label{eq24}&\Rightarrow& x_1^{\,\sharp\, (nv)}=\left(x_1^{\,\sharp\, n}\right)^{\,\sharp\, v}=0 \quad\mbox{for some$v\in \mathcal{Z}_{>0}$.}
\end{eqnarray}

It follows from (\ref{eq23}) and (\ref{eq24}) that $x\in nilrad ^{\sharp}(R)$. Hence, 
$x+nilrad ^{\sharp}(R)=nilrad ^{\sharp}(R)$. This proves (ii).

\hfill\raisebox{1mm}{\framebox[2mm]{}}

\bigskip
We now characterize the nil radical of a left commutative rng by using Hu-Liu prime ideals. 

\begin{proposition}\label{pr2.3} The nil radical of a left commutative rng $R$ is the intersection of the Hu-Liu prime ideals of $R$.
\end{proposition}

\medskip
\noindent
{\bf Proof} In this proof, $R$ denotes a left commutative rng with an left identity $1^{\ell}$, the $\varepsilon$-component $x_{\varepsilon}$ of $x\in R$ and the $\varepsilon$-part 
$I_{\varepsilon}$ of an ideal $I$ of $R$ are always induced by the left identity $1^{\ell}$.

\medskip
Let $x$ be any element of $nilrad ^{\sharp}(R)$. Then $x^m=0$ and $x_1^{\,\sharp\, n}=0$ for some $m$, $n\in\mathcal{Z}_{>0}$. Let $P$ be any Hu-Liu prime ideal of $R$. Since $x^m=0\in P+\hbar ^+(R)$, we have $x\in P+\hbar ^+(R)=P_0\oplus \hbar ^+(R)$ or
\begin{equation}\label{eq25} 
x_0\in P_0\subseteq P.
\end{equation}

If $P\not\supseteq \hbar ^+(R)$, then $\left(\, \hbar ^+\left(\displaystyle\frac{R}{P}\right) , \, + , \, \,\sharp\, \, \right)$ is a domain by Proposition~\ref{pr1.1}(v). Using $x_1^{\,\sharp\, n}=0$, we get
$$
(x_1+P)^{\,\sharp\, n}=x_1^{\,\sharp\, n}+P=0+P=P.
$$
Thus $x_1+P$ is a nilpotent element of the domain 
$\hbar ^+\left(\displaystyle\frac{R}{P}\right)$, which implies that $x_1+P$ is the zero element of the domain 
$\hbar ^+\left(\displaystyle\frac{R}{P}\right)$. Hence, we get
\begin{equation}\label{eq26} 
x_1\in P.
\end{equation}
If $P\supseteq \hbar ^+(R)$, then (\ref{eq26}) is obviously true.

It follows from (\ref{eq25}) and (\ref{eq26}) that $x=x_0+x_1\in P$. This proves that
\begin{equation}\label{eq27} 
nilrad ^{\sharp}(R)\subseteq \bigcap_{P\in spec^{\sharp}R} P.
\end{equation}

\medskip
Conversely, we prove that
\begin{equation}\label{eq28} 
z\ne nilrad ^{\sharp}(R) \,\Rightarrow\, z\ne \bigcap_{P\in spec^{\sharp}R} P.
\end{equation}

\medskip
\underline{\it Case 1:} $z^m\ne 0$ for all $m\in \mathcal{Z}_{>0}$, in which case, 
$z^m\ne \hbar ^+(R)$ for all $m\in \mathcal{Z}_{>0}$. Hence, $z+ \hbar ^+(R)$ is not a nilpotent element of the commutative ring $\displaystyle\frac{R}{\hbar ^+(R)}$; that is,
$$
z+ \hbar ^+(R)\not\in nilrad\left(\displaystyle\frac{R}{\hbar ^+(R)} \right)
=\bigcap _{\frac{I}{\hbar ^+(R)}\in spec \left( \frac{R}{\hbar ^+(R)}\right)} 
\left(\displaystyle\frac{I}{\hbar ^+(R)}\right).
$$
Hence, there exists a prime ideal $\displaystyle\frac{I}{\hbar ^+(R)}$ of the commutative ring $\displaystyle\frac{R}{\hbar ^+(R)}$ such that $z\not\in I$. Since $I$ is a Hu-Liu prime ideal by 
Proposition~\ref{pr1.1}(ii), (\ref{eq28}) holds in this case.

\medskip
\underline{\it Case 2:} $z_1^{\,\sharp\, n}\ne 0$ for all $n\in \mathcal{Z}_{>0}$, in which case, we consider the following set
$$
T:=\left\{\, J \, \left| \, \mbox{$J$ is an ideal of $R$ and $z_1^{\,\sharp\, n}\not\in J$ for all $n\in \mathcal{Z}_{>0}$}\right.\right\}.
$$
Since $\{0\}\in T$, $T$ is nonempty. Clearly, $(\,T, \, \subseteq \,)$ is a partially order set, where $\subseteq $ is the relation of set inclusion. If 
$\{\, J_{\lambda} \, |\, \lambda\in \Lambda \,\}$ is a nonempty totally ordered subset of $T$, then $\cup_{\lambda\in \Lambda} J_{\lambda}$ is an upper bound of 
$\{\, J_{\lambda} \, |\, \lambda\in \Lambda \,\}$ in $T$. By Zorn's Lemma, the partially ordered set $(\,T, \, \subseteq \,)$ has a maximal element $P$. We are going to prove that $P$ is a Hu-Liu prime ideal of $R$.

\medskip
Let $x$ and $y$ be two elements of $R$. First, if $x\not\in P+ \hbar ^+(R)$ and 
$y\not\in P+ \hbar ^+(R)$, then $x_0\not\in P_0$ and $y_0\not\in P_0$. Hence, we get
\begin{equation}\label{eq29}
P\subset P+x_0R\quad\mbox{and}\quad P\subset P+y_0R.
\end{equation}
By (\ref{eq1}), we have
$$
R(P+x_0R)\subseteq  P+x_0R\quad\mbox{and}\quad (P+x_0R)R\subseteq P+x_0R.
$$
Moreover, $(P+x_0R)\cap \hbar ^+(R)=P_1+x_0R_1$ is an ideal of 
$(\, \hbar ^+(R), \, + , \, \,\sharp\, \,)$. This proves that the subgroup 
$P+x_0R$ is an ideal of the left commutative rng $R$. Similarly, $P+y_0R$ is also an ideal of the $R$. It follows from (\ref{eq29}) that
$$
P+x_0R\not\in T\quad\mbox{and}\quad P+y_0R\not\in T,
$$
which imply that
$$
z_1^{\,\sharp\, u}\in P+x_0R\quad\mbox{and}\quad z_1^{\,\sharp\, v}\in P+y_0R
$$
or
$$
z_1^{\,\sharp\, u}\in P_1+x_0R_1\quad\mbox{and}\quad z_1^{\,\sharp\, v}\in P_1+y_0R_1\quad\mbox{for some $u$, $v\in \mathcal{Z}_{>0}$.}
$$
By (\ref{eq3}), we have
\begin{eqnarray*}
&&z_1^{\,\sharp\, (u+v)}=z_1^{\,\sharp\, u}\,\sharp\, z_1^{\,\sharp\, v}\in 
(P_1+x_0R_1)\,\sharp\, (P_1+y_0R_1)\\&&\\
&\subseteq& 
\underbrace{P_1\,\sharp\, P_1+ P_1\,\sharp\, (y_0R_1)+
(x_0R_1)\,\sharp\, P_1}_{\mbox{This is a subset of $P$}}+ (x_0R_1)\,\sharp\, (y_0R_1)\\&&\\
&\subseteq& P+x_0y_0(R_1\,\sharp\, R_1)\subseteq P+x_0y_0R,
\end{eqnarray*}
which implies that
\begin{equation}\label{eq30} 
(xy)_0=x_0y_0\not\in P.
\end{equation}
Hence, $xy\not\in P+\hbar ^+(R)$. This proves that
\begin{equation}\label{eq31} 
\mbox{$x\not\in P+ \hbar ^+(R)$ and $y\not\in P+ \hbar ^+(R)$} \,\Rightarrow\,
xy\not\in P+\hbar ^+(R).
\end{equation}

\medskip
Similarly, we can prove 
\begin{equation}\label{eq32} 
\mbox{$x\not\in P+ \hbar ^+(R)$ and $y\not\in P$} \,\Rightarrow\,
xy\not\in P.
\end{equation}
and
\begin{equation}\label{eq33} 
\mbox{$P_1=P\cap \hbar ^+(R)$ is a prime ideal of  
$(\, \hbar ^+(R),\, +,\, \,\sharp\, \,)$. }
\end{equation}

\medskip
By (\ref{eq31}), (\ref{eq32}) and (\ref{eq33}), $P$ is a Hu-Liu prime ideal. Since 
$z_1\not\in P$, (\ref{eq28}) also holds in Case 2.

\medskip
It follows from (\ref{eq27}) and (\ref{eq28}) that Proposition~\ref{pr2.3} is true.

\hfill\raisebox{1mm}{\framebox[2mm]{}}

\bigskip
If $I$ is an ideal of a left commutative rng $R$, then the {\bf radical} $\sqrt[\sharp]{I}$ of $I$ is defined by
$$
\sqrt[\sharp]{I}:=\{\, x\in R \, | \, \mbox{$x_0^m\in I\cap R1^{\ell}$ and $x_1^{\,\sharp\, n}\in I\cap\hbar ^+(R)$ for some $m$, $n\in\mathcal{Z}_{>0}$} \,\},
$$
where $x_\varepsilon$ is the $\varepsilon$-component of $x$ induced by the left identity $1^{\ell}$ of $R$. Since
$$
nilrad ^{\sharp}\left(\displaystyle\frac{R}{I}\right)=\displaystyle\frac{\sqrt[\sharp]{I}+I}{I},
$$
$\sqrt[\sharp]{I}$ is an ideal of $R$. An ideal $I$ of a left commutative rng $R$ is called a {\bf radical ideal} if $\sqrt[\sharp]{I}=I$.

\medskip
The next proposition is a corollary of Proposition~\ref{pr2.3}.

\begin{proposition}\label{pr2.4} If $I$ is an ideal of a left commutative rng $R$ and $I\ne R$, then
$$\sqrt[\sharp]{I}=\bigcap _{\mbox{$P\in spec^{\sharp}R$ and $P\supseteq I$}} P.$$
\end{proposition}

\medskip
\noindent
{\bf Proof} By Proposition~\ref{pr2.3}, we have
\begin{eqnarray*}
&&x\in \sqrt[\sharp]{I}\\
&\Leftrightarrow& x+I\in nilrad ^{\sharp}\left(\frac{R}{I}\right)
=\bigcap _{\mbox{$\frac{P}{I}\in spec^{\sharp}\left(\frac{R}{I}\right)$}}\frac{R}{I}\\
&\Leftrightarrow& x\in \bigcap _{\mbox{$P\in spec^{\sharp}R$ and $P\supseteq I$}} P.
\end{eqnarray*}

\hfill\raisebox{1mm}{\framebox[2mm]{}}

\medskip
\section{The Expansion of the Zarisky Topology}

Let $(\, R, \, + , \, \cdot , \, \,\sharp\, \,)$ be a left commutative rng. For an ideal $I$ of $R$, we define a subset $\mathcal{V}(I)$ of $spec^{\sharp}R$ by
\begin{equation}\label{eq34} 
\mathcal{V}(I): =\{\, P \, |\, \mbox{$P\in spec^{\sharp}R$ and $P\supseteq I$} \,\}.
\end{equation}

\begin{proposition}\label{pr3.1} Let $R$ be a left commutative rng $R$.
\begin{description}
\item[(i)] $\mathcal{V}(0)=spec^{\sharp}R$ and $\mathcal{V}(R)=\emptyset$.
\item[(ii)] $\mathcal{V}(I)\cup \mathcal{V}(J)=\mathcal{V}(I\cap J)$, where $I$ and $J$ are two ideals of $R$.
\item[(iii)] $\displaystyle\bigcap _{\lambda \in \Lambda} 
\mathcal{V}\Big(\stackrel{\lambda}{I}\Big)=
\mathcal{V}\Big(\displaystyle\sum _{\lambda \in \Lambda}\stackrel{\lambda}{I}\Big)$, where 
$\Big\{\left.\,\stackrel{\lambda}{I}\, \right| \, \lambda \in \Lambda\, \Big\}$ is a set of ideals of $R$.
\end{description}
\end{proposition}

\medskip
\noindent
{\bf Proof} Since (i) and (iii) are clear, we need only to prove (ii).

\medskip
By (\ref{eq34}), we have
\begin{equation}\label{eq35} 
\mathcal{V}(I)\cup \mathcal{V}(J)\subseteq \mathcal{V}(I\cap J).
\end{equation}

\medskip
Conversely, we prove that
\begin{equation}\label{eq36} 
\mbox{$P\supseteq I\cap J$ and $P\not\supseteq I$} \Rightarrow P\supseteq J.
\end{equation}

In the following proof, we assume that the even part and  odd part of an ideal are always induced by the same left identity. Since $P\not\supseteq I$, we have two cases.

\medskip
\underline{\it Case 1:} $P\not\supseteq I_0$, in which case, there exists $x_0\in I_0$ and $x_0\not\in P_0$. Hence. $x_0J\subseteq I\cap J\subseteq P$, which implies $P\supseteq J$ by (\ref{eq15}).

\medskip
\underline{\it Case 2:} $P\supseteq I_0$ and $P\not\supseteq I_1$, in which case, there exists $x_1\in I_1$ and $x_1\not\in P$. Hence. $J_0x_1\subseteq I\cap J\subseteq P$. This fact and (\ref{eq15}) imply that $J_0\subseteq P+\hbar ^+(R)$ or
\begin{equation}\label{eq37} 
J_0\subseteq P_0.
\end{equation}

Also, we have 
\begin{equation}\label{eq38} 
x_1\,\sharp\, J_1\subseteq (I\cap J)_1\subseteq P_1.
\end{equation}
Since $P\not\supseteq I_1$, $P\not\supseteq \hbar ^+(R)$. Thus, $P_1=P\cap \hbar ^+(R)$ is a prime ideal of the commutative ring $(\, \hbar ^+(R) , \, + , \, \,\sharp\, \, )$. Hence, (\ref{eq38}) implies that
\begin{equation}\label{eq39} 
J_1\subseteq P_1.
\end{equation}

By (\ref{eq37}) and (\ref{eq39}), we get $J=J_0+J_1\subseteq P_0+P_1=P$. This proves that (\ref{eq36}) is also true in Case 2.

\medskip
Using (\ref{eq36}), we have
\begin{equation}\label{eq40} 
\mathcal{V}(I)\cup \mathcal{V}(J)\supseteq \mathcal{V}(I\cap J).
\end{equation}

It follows from (\ref{eq35}) and (\ref{eq40}) that (ii) is true.

\hfill\raisebox{1mm}{\framebox[2mm]{}}

\bigskip
Let $R$ be a left commutative rng. By Proposition~\ref{pr3.1}, the collection
$$
\mathcal{V}: =\{\, \mathcal{V}(I) \, |\, \mbox{$I$ is an ideal of $R$} \, \}
$$
of subsets of $spec^{\sharp}R$ satisfies the axioms for closed sets in a topological space. The topology on $spec^{\sharp}R$ having the elements of $\mathcal{V}$ as closed sets is called the {\bf expanded Zariski topology}. The collection
$$
\mathcal{D}: =\{\, \mathcal{D}(I) \, |\, \mbox{$I$ is an ideal of $R$} \, \}
$$
consists of the open sets of the expanded Zariski topology on $spec^{\sharp}R$, where
$$
\mathcal{D}(I): =spec^{\sharp}R\setminus \mathcal{V}(I)=
\{\, P \, |\, \mbox{$P\in spec^{\sharp}R$ and $P\not\supseteq I$} \, \}.
$$
Since $spec^{\sharp}_0R=\mathcal{V}(\hbar ^+(R))$, the even spectrum $spec^{\sharp}_0R$ of $R$ is a closed subspace of $spec^{\sharp}R$ and the odd spectrum $spec^{\sharp}_1R$ of $spec^{\sharp}R$ is an open subspace of $spec^{\sharp}R$. 
We define a binary relation $\stackrel{\sharp}{\sim}$ on $spec^{\sharp}_1R$ by
$$
P\stackrel{\sharp}{\sim} Q \Leftrightarrow P\cap \hbar ^+(R)=Q\cap \hbar ^+(R)\quad\mbox{for $P$, $Q\in spec^{\sharp}_1R$}.
$$
It is clear that $\stackrel{\sharp}{\sim}$ is an equivalence relation. The equivalence class containing $P\in spec^{\sharp}_1R$ is denoted by $[P]$. The quotient topology 
with respect to the equivalence relation $\stackrel{\sharp}{\sim}$ is called the 
{\bf odd quotient topology} of  $spec^{\sharp}_1R$ and is denoted by $spec^{\sharp}_1R/\stackrel{\sharp}{\sim}$. 

\medskip
The next proposition gives the basic properties of the closed subspace of $spec^{\sharp}R$ and the odd quotient topology of  $spec^{\sharp}_1R$.

\begin{proposition}\label{pr3.2} Let $(\, R, \, + , \, \cdot , \, \,\sharp\, \,)$ be a left commutative rng $R$.
\begin{description}
\item[(i)] The map
$$
\phi_0^R: \displaystyle\frac{P}{\hbar ^+(R)}\mapsto P \quad
\mbox{for $\displaystyle\frac{P}{\hbar ^+(R)}\in spec\left(\displaystyle\frac{R}{\hbar ^+(R)}\right)$}
$$
is a homeomorphism from the spectrum of the commutative ring $\displaystyle\frac{R}{\hbar ^+(R)}$
onto the closed subspace $spec^{\sharp}_0R$ of $spec^{\sharp}_1R$.
\item[(ii)] The map
$$
\phi_1^R: [P]\mapsto P\cap \hbar ^+(R) \quad
\mbox{for $[P]\in spec^{\sharp}_1R/\stackrel{\sharp}{\sim}$}
$$
is a homeomorphism  from the odd quotient topology space $spec^{\sharp}_1R/\stackrel{\sharp}{\sim}$ onto the subspace 
$\{\, P\cap \hbar ^+(R) \, |\, P\in spec^{\sharp}_1R \,\}$ of the spectrum of the commutative ring $(\,\hbar ^+(R), \, +, \, \sharp \,)$.
\end{description}
\end{proposition}

\medskip
\noindent
{\bf Proof} Clear.

\hfill\raisebox{1mm}{\framebox[2mm]{}}

\bigskip
For any commutative ring $\bar{R}$, there exists a left commutative rng $R$ such that 
$\bar{R}\simeq\displaystyle\frac{R}{\hbar ^+(R)}$ as rings. It follows from this fact and Proposition\ref{pr3.2} that the points in the spectrum of a commutative ring appear as the points in the even spectrum of a left commutative rng, and the new points introduced to expand the Zariski topology  consist of the odd spectrum of the left commutative rng.

\begin{proposition}\label{pr3.3} Let $R$ and $S$ be left commutative rngs. If $f: R\to S$ is a left commutative rng homomorphism, then the map
$$
f^{\sharp}: P\mapsto f^{-1}(P) \quad\mbox{for $P\in spec^{\sharp}S$}
$$
is continuous from $spec^{\sharp}S$ to $spec^{\sharp}R$, and the restriction map 
$\left.f^{\sharp}\right|spec^{\sharp}_{\varepsilon}S$ is continuous from the subspace $spec^{\sharp}_{\varepsilon}S$ of $spec^{\sharp}S$ to the subspace $spec^{\sharp}_{\varepsilon}R$ of $spec^{\sharp}R$.
\end{proposition}

\medskip
\noindent
{\bf Proof} First, we prove
\begin{equation}\label{eq41} 
P\in spec^{\sharp}S\Rightarrow f^{-1}(P)\in spec^{\sharp}R.
\end{equation}

Let $1_R^{\ell}$ be a left identity of $R$ such that $1_S^{\ell}:=f(1_R^{\ell})$ is a left identity of $S$. In the following proof, the $\varepsilon$-component of an element of $R$ are always induced by $1_R^{\ell}$, and the $\varepsilon$-component of an element of $S$ are always induced by $1_S^{\ell}$. 

For $x_{\varepsilon}$, $y_{\varepsilon}\in R_{\varepsilon}$ with $\varepsilon=0$ and $1$, we have
\begin{eqnarray*}
&&x_0y_{\varepsilon}\in \big(f^{-1}(P)\big)_{\varepsilon}=f^{-1}(P)\cap R_{\varepsilon}\\
&\Rightarrow& f(x_0)f(y_{\varepsilon})=f(x_0y_{\varepsilon})\in P\cap f(R_{\varepsilon})
\subseteq P\cap S_{\varepsilon}=P_{\varepsilon}\\
&\Rightarrow& f(x_0)\in P_0\quad\mbox{or}\quad f(y_{\varepsilon})\in P_{\varepsilon}\\
&\Rightarrow& x_0\in f^{-1}(P)\cap R_0=\big(f^{-1}(P)\big)_0\quad\mbox{or}\quad 
y_{\varepsilon}\in f^{-1}(P)\cap R_{\varepsilon}=\big(f^{-1}(P)\big)_{\varepsilon}.
\end{eqnarray*}

This proves that 
\begin{equation}\label{eq42} 
x_0y_{\varepsilon}\in \big(f^{-1}(P)\big)_{\varepsilon}\Rightarrow 
x_0\in \big(f^{-1}(P)\big)_0 \quad\mbox{or}\quad y_{\varepsilon}\in 
\big(f^{-1}(P)\big)_{\varepsilon}.
\end{equation}

If $f^{-1}(P)\supseteq R_1$, then (\ref{eq41}) is true by (\ref{eq42}).

If $f^{-1}(P)\not\supseteq R_1$, then $P\not\supseteq S_1$; otherwise, $R_1\subseteq f^{-1}(S_1)
\subseteq f^{-1}(P)$. Hence, $P\cap S_1$ is a prime ideal of $S_1$. Thus, we have
\begin{eqnarray*}
&& x_1\,\sharp\, y_1\in f^{-1}(P)\cap R_1\\
&\Rightarrow& f(x_1)\,\sharp\, f(y_1)\in P\cap f(R_1)\subseteq P\cap S_1\\
&\Rightarrow& f(x_1)\in P\cap S_1\quad\mbox{or}\quad f(y_1)\in P\cap S_1\\
&\Rightarrow& x_1\in f^{-1}(P)\cap R_1\quad\mbox{or}\quad y_1\in f^{-1}(P),
\end{eqnarray*}
which proves that
\begin{equation}\label{eq43} 
\mbox{$f^{-1}(P)\not\supseteq R_1$ $\Rightarrow$ $f^{-1}(P)\cap R_1$ is a prime ideal of $R_1$.}
\end{equation}
By (\ref{eq42}) and (\ref{eq43}), (\ref{eq41}) is also true.

It follows from (\ref{eq41}) that the map $f^{\sharp}$ is well-defined.

\medskip
Next, we prove that
\begin{equation}\label{eq44} 
\mbox{$I$ is an ideal of $R$ $\Rightarrow$  $\big(f^{\sharp}\big)^{-1}\big( \mathcal{V}(I)\big)=
\mathcal{V}(<f(I)>)$,}
\end{equation}
where $<f(I)>$ is the ideal of $S$ generated by $f(I)$. Since
\begin{eqnarray*}
&&Q\in \big(f^{\sharp}\big)^{-1}\big( \mathcal{V}(I)\big)\\
&\Leftrightarrow &f^{-1}(Q)=f^{\sharp}(Q)\in \mathcal{V}(I)\\
&\Leftrightarrow& f^{-1}(Q)\supseteq I\\
&\Leftrightarrow &Q\supseteq f(I)\\
&\Leftrightarrow& Q\supseteq <f(I)>\\
&\Leftrightarrow& Q\in \mathcal{V}(<f(I)>),
\end{eqnarray*}
we know that (\ref{eq44}) is true.

\medskip
By (\ref{eq44}), the pre-image of any closed set in $spec^{\sharp}R$ under $f^{\sharp}$ is a closed set in $spec^{\sharp}S$. Hence, the map $f^{\sharp}$ is continuous.

\medskip
Finally, let $1^{\sharp}_R$ and $1^{\sharp}_S$ be the local identity of $R$ and $S$, respectively. If $I$ is an ideal of $R$, then
$$
f\big(I+\hbar ^+(R)\big)=f(I)+f\big(\hbar ^+(R)\big)\subseteq f(I)+\hbar ^+(S)
\subseteq <f(I)>+\hbar ^+(S)
$$
which implies that
\begin{equation}\label{eq45}
<f\big(I+\hbar ^+(R)\big)>\subseteq <f(I)>+\hbar ^+(S).
\end{equation}

Conversely, since $1^{\sharp}_S=f(1^{\sharp}_R)\in f\big(\hbar ^+(R)\big)\subseteq 
f\big(I+\hbar ^+(R)\big)$, we get $\hbar ^+(S)$ $\subseteq <f\big(I+\hbar ^+(R)\big)>$. It follows that
\begin{equation}\label{eq46}
<f\big(I+\hbar ^+(R)\big)>\supseteq <f(I)>+\hbar ^+(S).
\end{equation}

By (\ref{eq45}) and (\ref{eq46}), we get
\begin{equation}\label{eq47}
<f\big(I+\hbar ^+(R)\big)>= <f(I)>+\hbar ^+(S).
\end{equation}

It follows from (\ref{eq44}) and (\ref{eq47}) that
$$
\big(f^{\sharp}\big)^{-1}\big( \mathcal{V}(I+\hbar ^+(R))\big)=
\mathcal{V}\big(<f(I)>+\hbar ^+(S)\big),
$$
which proves that the restriction map 
$\left.f^{\sharp}\right|spec^{\sharp}_{\varepsilon}S$ is continuous.

\hfill\raisebox{1mm}{\framebox[2mm]{}}

\bigskip
By Proposition~\ref{pr3.3}, the pair of maps $R\mapsto spec^{\sharp}R$, $f\mapsto f^{\sharp}$ define a contravariant functor from the category of left commutative rngs (left commutative rng homomorphisms as morphisms) into the category of topological spaces (continuous maps as morphisms). Moreover, if $f: R\to S$ is a left commutative rng homomorphism, then $f$ induces a commutative ring homomorphism $\bar{f}: \displaystyle\frac{R}{\hbar ^+(R)}\to 
\displaystyle\frac{S}{\hbar ^+(S)}$, and the following two diagrams
$$\begin{array}{ccc}
R&\stackrel{f}{\longrightarrow}&S\\
\shortstack{}\Bigg\downarrow&&\Bigg\downarrow\shortstack{}\\&&\\
\displaystyle\frac{R}{\hbar ^+(R)}&\stackrel{\bar{f}}{\longrightarrow}&\displaystyle\frac{S}{\hbar ^+(S)}
\end{array}\quad\mbox{and}\quad 
\begin{array}{ccc}
spec^{\sharp}S&\stackrel{f^{\sharp}}{\longrightarrow}&spec^{\sharp}R\\
\shortstack{$\phi_0^S$}\Bigg\uparrow&&\Bigg\uparrow\shortstack{$\phi_0^R$}\\&&\\
spec\left(\displaystyle\frac{S}{\hbar ^+(S)}\right)&\stackrel{\bar{f}^\ast}{\longrightarrow}&spec\left(\displaystyle\frac{R}{\hbar ^+(R)}\right)
\end{array}
$$
are commutative.

\medskip
Other properties of left commutative rngs can be found in \cite{Liu5}.

\end{document}